# A geometric approach to the linear trace Harnack inequality for the Ricci flow *

Bennett Chow and Sun-Chin Chu
School of Mathematics
University of Minnesota
Minneapolis, MN 55455†

## 1 Introduction

In [H], Richard Hamilton proved a matrix Harnack inequality for the Ricci flow, a consequence of which is the following trace Harnack inequality

**Theorem A (Hamilton).** *If $(M, g)$ is a complete solution to the Ricci flow*

$$\frac{\partial}{\partial t}g_{ij} = -2R_{ij} \qquad (1.1)$$

*with nonnegative curvature operator and bounded curvature, then for any 1-form $V$*

$$\frac{\partial}{\partial t}R + \frac{R}{t} + 2\nabla R \cdot V + 2R_{ij}V^i V^j \geq 0. \qquad (1.2)$$

*In particular, taking $V = 0$, we have*

$$\frac{\partial}{\partial t}(tR) \geq 0.$$

This trace inequality turns out to be a special case of the following linear Harnack inequality, which was later proved by the first author and Hamilton [CH].

**Theorem B.** *Under the same hypotheses as* Theorem A, *if $h$ is a nonnegative symmetric 2-tensor satisfying*

$$\frac{\partial}{\partial t}h_{ij} = \Delta h_{ij} + 2R_{pijq}h_{pq} - R_{iq}h_{jq} - R_{jq}h_{iq}, \qquad (1.3)$$

---

where the Laplacian and curvature are with respect to the metric evolving under the Ricci flow, then

$$Z := div(div(h)) + Rc \cdot h + 2div(h) \cdot V + h_{ij}V^i V^j + \frac{H}{2t} \geq 0, \qquad (1.4)$$

where $H = g^{ij}h_{ij}$. In particular, taking $h = Rc$, we obtain Theorem A as a special case.

On the other hand, in [CC], the authors showed that Hamilton's matrix Harnack quadratic for the Ricci flow is actually the Riemann curvature tensor of a connection on the space-time manifold $M \times [0, T)$ compatible with a degenerate metric on space-time. In particular, recall the space-time metric and connection defined in [CC]. The degenerate metric $\tilde{g}$ on the cotangent bundle $T^*\tilde{M}$ is defined by

$$\tilde{g}^{ij} = \begin{cases} g^{ij} & \text{if } i, j \geq 1 \\ 0 & \text{if } i = 0 \text{ or } j = 0. \end{cases}$$

Associated to this metric is the space-time connection $\tilde{\nabla}$ defined by

$$\tilde{\Gamma}_{ij}^k = \begin{cases} \Gamma_{ij}^k & \text{if } i, j, k \geq 1 \\ 0 & \text{if } k = 0 \text{ and } i, j \geq 0 \\ -R_i^j & \text{if } i = 0 \text{ and } j, k \geq 1 \\ -\frac{1}{2}\nabla^k R & \text{if } i = j = 0 \text{ and } k \geq 1. \end{cases}$$

This connection is compatible with the metric in the usual sense that

$$\tilde{\nabla}\tilde{g} = 0.$$

Moreover, it has the special property that

$$\frac{\partial}{\partial t}\tilde{\Gamma}_{ij}^k = -\tilde{g}^{kl}(\tilde{\nabla}_i \tilde{R}_{jl} + \tilde{\nabla}_j \tilde{R}_{il} - \tilde{\nabla}_l \tilde{R}_{ij}),$$

which is formally the same equation as that satisfied by the Levi-Civita connection $\nabla$ of the metric $g$ evolving under the Ricci flow. In [CC] it was shown that the Riemann curvature tensor of the connection $\tilde{\nabla}$ is the same as Hamilton's matrix Harnack quadratic. Similarly, the Ricci tensor of the $\tilde{\nabla}$ is the same as the trace Harnack quadratic.

**Theorem C.** *Given a 1-form $W_i$ and a 2-form $U_{ij}$, let*

$$Q = \left[\Delta R_{ij} - \frac{1}{2}\nabla_i\nabla_j R + 2R_{kijl}R_{kl} - R_{ik}R_{kj}\right]W_i W_j - 2(\nabla_i R_{jk} - \nabla_j R_{ik})U_{ij}W_k$$
$$+ R_{ijkl}U_{ij}U_{lk}$$

*denote Hamilton's Harnack quadratic. We have*

$$Q = \tilde{g}^{ip}\tilde{R}^l_{pjk}T^j_i T^k_l,$$

*where*

$$T^j_i = \begin{cases} U^j_i & \text{if } i,j \geq 1 \\ W_i & \text{if } j = 0 \\ 0 & \text{if } i = 0. \end{cases}$$

In this paper, we shall show that one can approach the linear trace Harnack inequality from this point of view. In section 2, we observe that the linear trace Harnack quadratic $Z$ given by (1.4) is equivalent to $\tilde{h}$, the natural extension to space-time of the symmetric 2-tensor $h$ given by (1.3). We then recall the space-time formulation in [CC] and show that the equation for $Z$ (i.e., $\tilde{h}$) derived in [CH] is the heat equation, using the Lichnerowicz Laplacian, in the space-time formulation. That is, the evolution equation for $\tilde{h}$ is the exact space-time analogue of the evolution equation for $h$. In section 3, we generalize the results of [CC] to the case of the Ricci flow modified by an arbitrary one-parameter family of diffeomorphism. In particular we define a suitable space-time connection and show that it satisfies the modified Ricci flow for degenerate metrics. The space-time formulation needs to be done in this generality in order to linearize the Ricci flow using DeTurck's trick. In section 4, we linearize the Ricci flow by considering a one-parameter family of Ricci flows modified by DeTurck's trick in the variation direction. We show that the variation of the metric satisfies (1.3) and the variation of the space-time connection satisfies formally the same equation as the Levi-Civita connection of the space metric, where $h$ and $\nabla$ are replaced by $\tilde{h}$ and $\tilde{\nabla}$. Both of these equations rely on using DeTurck's trick. In section 5, we extend DeTurck's trick to the space-time connection and show that the variation of the modified space-time Ricci tensor is given by the space-time Lichnerowicz Laplacian of $\tilde{h}$. This supports the viewpoint that $\tilde{h}$ is the variation of the pair $(\tilde{g}, \tilde{\nabla})$. Finally, in section 6, we show how the tensor $\tilde{h}$ and its evolution arises from taking the limit of Riemannian metrics on space-time.

## 2 The linear trace Harnack inequality

In this section we recall the computations in [CH] for the evolution of the linear trace Harnack quadratic

$$Z = div(div(\,h\,)) + Rc \cdot h + 2div(h) \cdot V + h_{ij}V^i V^j, \qquad (2.1)$$

and interpret them in terms of the connection and curvature tensor on space-time defined in [CC]. In sections 3 and 4, we shall explain why this interpretation holds.

From the computations in section 6 of [CH], we have

**Lemma 2.1** *Under the Ricci flow* (1.1) *and equation* (1.3) *for $h$,*

$$\frac{\partial}{\partial t}div(h)_i = \Delta div(h)_i + 2h_{pq}\nabla_i R_{pq} - 2h_{pq}\nabla_p R_{qi}$$
$$+ 2R_{pq}\nabla_p h_{qi} - R_{qi}div(h)_q$$
$$\frac{\partial}{\partial t}\left[div(div(h)) + Rc \cdot h\right] = \Delta\left[div(div(h)) + Rc \cdot h\right] + 4R_{ij}\nabla_j\left(div(h)_i\right)$$
$$+ 2h_{pq}\left(\Delta R_{pq} - \frac{1}{2}\nabla_p\nabla_q R + 2R_{pijq}R_{ij}\right)$$

We now interpret these computations in terms of the space-time formalism for the Ricci flow. Let $\tilde{M} = M \times [0,T)$ be the space-time manifold. Define a symmetric 2-tensor on $\tilde{M}$

$$\tilde{h} = \sum_{i,j=0}^{n} \tilde{h}_{ij} dx^i \otimes dx^j$$

by

$$\tilde{h}_{ij} = \begin{cases} h_{ij} & \text{if } i,j \geq 1 \\ div(h)_j & \text{if } i = 0 \\ div(div(h)) + Rc \cdot h & \text{if } i = j = 0, \end{cases}$$

where $\{x^i\}_{i=1}^{n}$ are local coordinates on $M$ and $x^0 = t$ is the time coordinate. We may now rewrite the linear trace Harnack quadratic (2.1) as

$$Z = \sum_{i,j=0}^{n} \tilde{h}_{ij} \bar{V}^i \bar{V}^j,$$

where $\bar{V} = V \oplus \frac{\partial}{\partial t}$. The conclusion of Theorem B may now be restated as

$$\sum_{i,j=0}^{n} \tilde{h}_{ij}(\tilde{V}^i \tilde{V}^j + \frac{1}{2t}\tilde{g}^{ij}) \geq 0,$$

for any space-time vector field $\tilde{V} = V \oplus \tilde{V}^0 \frac{\partial}{\partial t}$.

Using the space-time connection, we can state the following special property of the symmetric 2-tensor $\tilde{h}$.

**Lemma 2.2** *For all $0 \leq i \leq n$, we have*

$$\tilde{h}_{i0} = \tilde{g}^{jk}\tilde{\nabla}_j \tilde{h}_{ki}. \tag{2.2}$$

PROOF. We only verify the case $i = 0$; the other case where $i \geq 1$ is even easier. Using the definitions of $\tilde{\nabla}$ and $\tilde{h}$, we compute

$$\begin{aligned}
\tilde{g}^{jk}\tilde{\nabla}_j \tilde{h}_{k0} &= g^{jk}\nabla_j div(h)_k - g^{jk}\tilde{\Gamma}^p_{j0}\tilde{h}_{kp} \\
&= div(div(h)) + g^{jk}R^p_j h_{kp} \\
&= \tilde{h}_{00}.
\end{aligned}$$

**Remark.** Formula (2.2) is analogous to the formulas

$$\tilde{R}^l_{ij0} = \tilde{g}^{pq}\tilde{\nabla}_p \tilde{R}^l_{ijq}$$

for $i, j, l \geq 0$, and

$$\tilde{R}_{i0} = \tilde{g}^{jk}\tilde{\nabla}_j \tilde{R}_{ki}$$

for $i \geq 0$, which were proved in [CC]. One may consider these formulas as defining the extension of a space tensor to the corresponding space-time tensor. In particular, we could have used Lemma 2.2 to define $\tilde{h}$ extending $h$.

The starting point for the space-time approach to the linear trace Harnack inequality is the following observation.

**Proposition 2.3** *Under the Ricci flow, if $h_{ij}$ satisfies*

$$\frac{\partial}{\partial t}h_{ij} = \Delta h_{ij} + 2R_{pijq}h_{pq} - R_{iq}h_{jq} - R_{jq}h_{iq}, \tag{2.3}$$

*then the associated space-time symmetric 2-tensor $\tilde{h}_{ij}$ satisfies*

$$\tilde{\nabla}_0 \tilde{h}_{ij} = \tilde{\Delta}\tilde{h}_{ij} + 2\tilde{h}_{pq}\tilde{R}^q_{pij}, \tag{2.4}$$

*or equivalently,*

$$\frac{\partial}{\partial t}\tilde{h}_{ij} = \tilde{\Delta}\tilde{h}_{ij} + 2\tilde{h}_{pq}\tilde{R}^q_{pij} - \tilde{R}_{iq}\tilde{h}_{jq} - \tilde{R}_{jq}\tilde{h}_{iq}, \tag{2.5}$$

*for all $i, j \geq 0$, where $\tilde{\Delta} = \tilde{g}^{ij}\tilde{\nabla}_i\tilde{\nabla}_j$ and $\tilde{R}^q_{pij}$ is the Riemann curvature (3,1)-tensor of $\tilde{\nabla}$.*

PROOF. When $i, j \geq 1$, equation (2.4) follows from (2.3) and the formula

$$\tilde{\nabla}_0 \tilde{h}_{ij} = \partial_0 \tilde{h}_{ij} - \tilde{\Gamma}^p_{0i} \tilde{h}_{pj} - \tilde{\Gamma}^p_{j0} \tilde{h}_{ip}$$
$$= \frac{\partial}{\partial t} h_{ij} + R^p_i h_{pj} + R^p_j h_{ip}.$$

For $i \geq 1$ and $j = 0$, (2.4) follows from Lemma 2.1, (i) and the equations

$$\tilde{\nabla}_0 \tilde{h}_{0i} = \partial_0 \tilde{h}_{0i} - \tilde{\Gamma}^p_{00} \tilde{h}_{pi} - \tilde{\Gamma}^p_{0i} \tilde{h}_{0p}$$
$$= \frac{\partial}{\partial t} div(h)_i + \frac{1}{2} \nabla^p R h_{pi} + R^p_i div(h)_p,$$

$$\tilde{\Delta} \tilde{h}_{0i} = \sum_{p=1}^n \tilde{\nabla}_p \tilde{\nabla}_p \tilde{h}_{0i}$$
$$= \nabla_p \left[ \nabla_p div(h)_i + R^q_p h_{qi} \right] - \tilde{\Gamma}^q_{p0} \tilde{h}_{qi}$$
$$= \Delta div(h)_i + \frac{1}{2} \nabla^p R h_{pi} + 2 R_{pq} \nabla_p h_{qi},$$

and

$$\tilde{R}^q_{p0i} = \nabla_i R^q_p - \nabla^q R_{pi}$$

(see [CC], Theorem 2.2, (B4) for the last equation.) Finally, for $i = j = 0$, equation (2.4) follows from combining Lemma 2.1, (ii) with

$$\tilde{\nabla}_0 \tilde{h}_{00} = \frac{\partial}{\partial t} \left[ div(div(h)) + Rc \cdot h \right] + \nabla^p R \cdot div(h)_p$$
$$\tilde{\Delta} \tilde{h}_{00} = \sum_{p=1}^n \tilde{\nabla}_p \tilde{\nabla}_p \tilde{h}_{00}$$
$$= \nabla_p \left[ \nabla_p (div(div(h)) + Rc \cdot h) + 2 R^q_p div(h)_q \right] - 2 \tilde{\Gamma}^q_{p0} \tilde{\nabla}_p \tilde{h}_{q0}$$
$$= \Delta \left[ div(div(h)) + Rc \cdot h \right] + 4 R_{pq} \nabla_p div(h)_q + 2 R_{pr} R_{rq} h_{pq} + \nabla^p R div(h)_p$$

and

$$\tilde{R}^q_{p00} = \Delta R_{pq} - \frac{1}{2} \nabla_p \nabla_q R + 2 R_{pijq} R_{ij} - R_{pr} R_{rq}$$

(see [CC], Theorem 2.2, (B5) for the last equation.) The proof of the proposition is complete.

# 3 Space-time formulation of the modified Ricci flow

In this section we extend the results of [CC] to the case of the Ricci flow modified by an arbitrary one-parameter family of diffeomorphisms. The space-time formulation needs to be done in this generality in order to linearize the Ricci flow using DeTurck's trick (see [D].) We consider the equation

$$\frac{\partial}{\partial t} g_{ij} = -2R_{ij} + \nabla_i V_j + \nabla_j V_i, \tag{3.1}$$

where $V = V(t)$ is an arbitrary one-parameter family of 1-forms. Recall that $\tilde{M} = M \times [0, T)$ and the space-time metric $\tilde{g}$ on the cotangent bundle $T^*\tilde{M}$ is the degenerate metric given by

$$\tilde{g}^{ij} = \begin{cases} g^{ij} \text{ if } i, j \geq 1 \\ 0 \text{ if } i = 0 \text{ or } j = 0, \end{cases} \tag{3.2}$$

where $x^0 = t$ is the time coordinate. We now define the space-time connection $\tilde{\nabla}$ by

(a) $\tilde{\Gamma}^k_{ij} = \Gamma^k_{ij}$ if $i, j, k \geq 1$
(b) $\tilde{\Gamma}^0_{ij} = 0$ if $i, j \geq 0$
(c) $\tilde{\Gamma}^k_{i0} = -R^k_i + \nabla_i V^k$ if $i, k \geq 1$
(d) $\tilde{\Gamma}^k_{00} = -\frac{1}{2}\nabla^k(R + |V|^2) + g^{kp}\frac{\partial}{\partial t}V_p$ if $k \geq 1$.

This definition, which agrees with the definition given in section 1 when $V = 0$, is natural for the following reasons. First, $\tilde{\nabla}$ is compatible with the metric $\tilde{g}$

$$\tilde{\nabla}_i \tilde{g}^{jk} = \frac{\partial}{\partial x^i}\tilde{g}^{jk} + \tilde{\Gamma}^j_{ip}\tilde{g}^{pk} + \tilde{\Gamma}^k_{ip}\tilde{g}^{jp} = 0$$

for all $i, j, k \geq 0$.

Second, let

$$\tilde{R}^l_{ijk} = \partial_i \tilde{\Gamma}^l_{jk} - \partial_j \tilde{\Gamma}^l_{ik} + \tilde{\Gamma}^p_{jk}\tilde{\Gamma}^l_{ip} - \tilde{\Gamma}^p_{ik}\tilde{\Gamma}^l_{jp}$$

and $\tilde{R}_{jk} = \tilde{R}^p_{pjk}$ denote the Riemann curvature and Ricci tensors of the space-time connection $\tilde{\nabla}$. Furthermore, extend the 1-form $V$ to space-time arbitrarily

$$\tilde{V} = V + h \cdot dt$$

($h$ is an arbitrary function,) i.e., $\tilde{V}_0 = h$. We then have **Theorem 3.1** *If $g$ satisfies the modified Ricci flow (3.1), then the space-time metric and connection satisfy the system*

$$\frac{\partial}{\partial t}\tilde{g}^{ij} = \tilde{g}^{ik}\tilde{g}^{jl}(2\tilde{R}_{kl} - \tilde{\nabla}_k\tilde{V}_l - \tilde{\nabla}_l\tilde{V}_k) \tag{3.3}$$

$$\frac{\partial}{\partial t}\tilde{\Gamma}^k_{ij} = -\tilde{\nabla}_i(\tilde{R}^k_j - \frac{1}{2}\tilde{\nabla}_j\tilde{V}^k - \frac{1}{2}\tilde{\nabla}^k\tilde{V}_j) - \tilde{\nabla}_j(\tilde{R}^k_i - \frac{1}{2}\tilde{\nabla}_i\tilde{V}^k - \frac{1}{2}\tilde{\nabla}^k\tilde{V}_i)$$
$$+\tilde{\nabla}^k(\tilde{R}_{ij} - \frac{1}{2}\tilde{\nabla}_i\tilde{V}_j - \frac{1}{2}\tilde{\nabla}_j\tilde{V}_i). \tag{3.4}$$

**Remark.** In [CC] we proved this result when $V = df$ is an exact 1-form and $\tilde{V}_0 = \frac{\partial}{\partial t}f$. There we conjectured that the result would hold when $V$ is an closed 1-form whose cohomology class is independent of time and where $\frac{\partial}{\partial t}V = d\tilde{V}_0$. It turns out, as the theorem says, we do not need to make any assumption on $V$ and $\tilde{V}_0$ may be taken arbitrarily.

The proof of the theorem relies on the following computations.

**Lemma 3.2**

$$\tilde{R}_{j0} - \frac{1}{2}\tilde{\nabla}_j\tilde{V}_0 - \frac{1}{2}\tilde{\nabla}_0\tilde{V}_j = \frac{1}{2}\nabla_j(R + |V|^2 - h) - \frac{1}{2}\frac{\partial}{\partial t}V_j \tag{3.5}$$

$$\tilde{R}_{00} - \tilde{\nabla}_0\tilde{V}_0 = \frac{1}{2}\frac{\partial}{\partial t}(R + |V|^2 - 2h) \tag{3.6}$$

We first show that the lemma implies the theorem.

**Proof of Theorem 3.1.** Using (3.1) and (3.2), it is easy to see that (3.3) holds for all $i, j \geq 0$. As for (3.4), the only nontrivial cases are
Case 1: $j = 0$ and $i, k \geq 1$
Case 2: $i, j = 0$ and $k \geq 1$.

For case 1, we compute

$$\frac{\partial}{\partial t}\tilde{\Gamma}^k_{i0} = \frac{\partial}{\partial t}(-R^k_i + \nabla_i V^k).$$

On the other hand, the RHS of (3.4) is given by

$$\text{RHS1} = -\tilde{\nabla}_i(\tilde{R}_0^k - \frac{1}{2}\tilde{\nabla}_0\tilde{V}^k - \frac{1}{2}\tilde{\nabla}^k\tilde{V}_0) - \tilde{\nabla}_0(\tilde{R}_i^k - \frac{1}{2}\tilde{\nabla}_i\tilde{V}^k - \frac{1}{2}\tilde{\nabla}^k\tilde{V}_i)$$
$$+ \tilde{\nabla}^k(\tilde{R}_{i0} - \frac{1}{2}\tilde{\nabla}_i\tilde{V}_0 - \frac{1}{2}\tilde{\nabla}_0\tilde{V}_i)$$
$$= -\nabla_i\left(\boxed{\frac{1}{2}\nabla^k(R+|V|^2-h)} - \frac{1}{2}g^{kp}\frac{\partial}{\partial t}V_p\right) + \tilde{\Gamma}_{i0}^p(R_p^k - \frac{1}{2}\nabla_p V^k - \frac{1}{2}\nabla^k V_p)$$
$$- \frac{\partial}{\partial t}(R_i^k - \frac{1}{2}\nabla_i V^k - \frac{1}{2}\nabla^k V_i) + \tilde{\Gamma}_{0i}^p(R_p^k - \frac{1}{2}\nabla_p V^k - \frac{1}{2}\nabla^k V_p)$$
$$- \tilde{\Gamma}_{0p}^k(R_i^p - \frac{1}{2}\nabla_i V^p - \frac{1}{2}\nabla^p V_i) + \nabla^k\left(\boxed{\frac{1}{2}\nabla_i(R+|V|^2-h)} - \frac{1}{2}\frac{\partial}{\partial t}V_i\right)$$
$$- g^{kl}\tilde{\Gamma}_{l0}^p(R_{ip} - \frac{1}{2}\nabla_i V_p - \frac{1}{2}\nabla_p V_i),$$

where we used Lemma 3.2, and the boxed terms cancel. Combining terms and using definition (c) for $\tilde{\Gamma}_{i0}^k$, we obtain

$$\text{RHS1} = \frac{1}{2}\frac{\partial}{\partial t}(g^{kp}\nabla_i V_p) - \frac{1}{2}(2R_{kp} - \nabla_k V_p - \nabla_p V_k)\nabla_i V_p + \frac{1}{2}g^{kp}(\frac{\partial}{\partial t}\Gamma_{ip}^q)V_q$$
$$- \frac{1}{2}\frac{\partial}{\partial t}(\nabla^k V_i) + \frac{1}{2}(2R_{kp} - \nabla_k V_p - \nabla_p V_k)\nabla_p V_i - \frac{1}{2}g^{kp}(\frac{\partial}{\partial t}\Gamma_{pi}^q)V_q$$
$$- \frac{\partial}{\partial t}(R_i^k - \frac{1}{2}\nabla_i V^k - \frac{1}{2}\nabla^k V_i) + (-R_i^p + \nabla_i V^p)(2R_{kp} - \nabla_k V_p - \nabla_p V_k)$$
$$- (-2R_{kp} + \nabla_k V_p + \nabla_p V_k)(R_{ip} - \frac{1}{2}\nabla_i V_p - \frac{1}{2}\nabla_p V_i)$$
$$= -\frac{\partial}{\partial t}(R_i^k - \nabla_i V^k) = \frac{\partial}{\partial t}\tilde{\Gamma}_{i0}^k,$$

since all of the rest of the terms cancel. This verifies (3.4) in case 1.

For case 2, we compute

$$\frac{\partial}{\partial t}\tilde{\Gamma}_{00}^k = \frac{\partial}{\partial t}\left(-\frac{1}{2}\nabla^k(R+|V|^2) + g^{kp}\frac{\partial}{\partial t}V_p\right).$$

The RHS of (3.4) is given by

$$\text{RHS2} = -2\tilde{\nabla}_0(\tilde{R}_0^k - \frac{1}{2}\tilde{\nabla}_0\tilde{V}^k - \frac{1}{2}\tilde{\nabla}^k\tilde{V}_0) + \tilde{\nabla}^k(\tilde{R}_{00} - \tilde{\nabla}_0\tilde{V}_0)$$
$$= -\frac{\partial}{\partial t}\left(\nabla^k(R+|V|^2-h) - g^{kp}\frac{\partial}{\partial t}V_p\right) + 2\tilde{\Gamma}_{00}^p(R_p^k - \frac{1}{2}\nabla_p V^k - \frac{1}{2}\nabla^k V_p)$$
$$- 2\tilde{\Gamma}_{0p}^k(\tilde{R}_0^p - \frac{1}{2}\tilde{\nabla}_0\tilde{V}^p - \frac{1}{2}\tilde{\nabla}^p\tilde{V}_0) + \frac{1}{2}\nabla^k\left(\frac{\partial}{\partial t}(R+|V|^2-2h)\right)$$
$$- g^{kl}\tilde{\Gamma}_{l0}^p(2\tilde{R}_{0p} - \tilde{\nabla}_0\tilde{V}_p - \tilde{\nabla}_p\tilde{V}_0).$$

Combining terms yields

$$
\begin{aligned}
\text{RHS2} &= \tfrac{\partial}{\partial t}\left[-\tfrac{1}{2}\nabla^k(R+|V|^2) + g^{kp}\tfrac{\partial}{\partial t}V_p\right] \\
&\quad - \tfrac{1}{2}(2R_{kp} - \nabla_k V_p - \nabla_p V_k)\nabla_p(R+|V|^2 - 2h) \\
&\quad + \left[-\tfrac{1}{2}\nabla^p(R+|V|^2) + g^{pq}\tfrac{\partial}{\partial t}V_q\right](2R_{kp} - \nabla_k V_p - \nabla_p V_k) \\
&\quad - (-2R_{kp} + \nabla_k V_p + \nabla_p V_k)(\nabla_p(R+|V|^2 - h) - \tfrac{\partial}{\partial t}V_p) \\
&= \tfrac{\partial}{\partial t}\left[-\tfrac{1}{2}\nabla^k(R+|V|^2) + g^{kp}\tfrac{\partial}{\partial t}V_p\right] = \tfrac{\partial}{\partial t}\tilde{\Gamma}^p_{00},
\end{aligned}
$$

which proves case 2 of (3.4), and hence the theorem.

We now return to the

**Proof of Lemma 3.2** First we compute

$$
\begin{aligned}
\tilde{R}^l_{ij0} &= \partial_i \tilde{\Gamma}^l_{j0} - \partial_j \tilde{\Gamma}^l_{i0} + \tilde{\Gamma}^p_{j0}\tilde{\Gamma}^l_{ip} - \tilde{\Gamma}^p_{i0}\tilde{\Gamma}^l_{jp} \\
&= \nabla_i(-R^l_j + \nabla_j V^l) - \nabla_j(-R^l_i + \nabla_i V^l) \\
&= -\nabla_i R^l_j + \nabla_j R^l_i + R^l_{ijp}V^p.
\end{aligned}
$$

Tracing implies

$$
\tilde{R}_{j0} = \tilde{R}^q_{qj0} = \tfrac{1}{2}\nabla_j R + R_{jp}V^p,
$$

where we used the contracted second Bianchi identity.

We also compute

$$
\begin{aligned}
\tilde{\nabla}_j \tilde{V}_0 + \tilde{\nabla}_0 \tilde{V}_j &= \nabla_j h - \tilde{\Gamma}^p_{j0}V_p + \tfrac{\partial}{\partial t}V_j - \tilde{\Gamma}^p_{0j}V_p \\
&= \nabla_j h + \tfrac{\partial}{\partial t}V_j + 2R^p_j V_p - \nabla_j |V|^2.
\end{aligned}
$$

Combining the two equations above,

$$
\tilde{R}_{j0} - \tfrac{1}{2}\tilde{\nabla}_j \tilde{V}_0 - \tfrac{1}{2}\tilde{\nabla}_0 \tilde{V}_j = \tfrac{1}{2}\nabla_j R - \tfrac{1}{2}\nabla_j h - \tfrac{1}{2}\tfrac{\partial}{\partial t}V_j + \tfrac{1}{2}\nabla_j |V|^2,
$$

which is equivalent to formula 1.

Second, we compute

$$
\begin{aligned}
\tilde{R}_{00} &= \tilde{R}^p_{p00} = \partial_p \tilde{\Gamma}^p_{00} - \partial_0 \tilde{\Gamma}^p_{p0} + \tilde{\Gamma}^q_{00}\tilde{\Gamma}^p_{pq} - \tilde{\Gamma}^q_{p0}\tilde{\Gamma}^p_{0q} \\
&= -\tfrac{1}{2}\Delta(R+|V|^2) + \nabla^p(\tfrac{\partial}{\partial t}V_p) - \tfrac{\partial}{\partial t}(-R + \nabla^p V_p) \\
&\quad - (-R^q_p + \nabla_p V^q)(-R^p_q + \nabla_q V^p) \\
&= -\tfrac{1}{2}\Delta R + \tfrac{\partial}{\partial t}R + R_{pq}V^p V^q - |Rc|^2,
\end{aligned}
$$

where to obtain the last equality we used the formula

$$-\frac{\partial}{\partial t}(\nabla^p V_p) = -\nabla^p(\frac{\partial}{\partial t}V_p) - (2R_{pq} - \nabla_p V_q - \nabla_q V_p)\nabla_p V_q + g^{pq}(\frac{\partial}{\partial t}\Gamma^r_{pq})V_r.$$

We also compute

$$\begin{aligned}\tilde{\nabla}_0 \tilde{V}_0 &= \frac{\partial}{\partial t}h - \tilde{\Gamma}^p_{00} V_p \\ &= \frac{\partial}{\partial t}h + \frac{1}{2}\nabla^p(R + |V|^2)V_p - g^{pq}\frac{\partial}{\partial t}V_q \cdot V_p \\ &= \frac{\partial}{\partial t}h + \frac{1}{2}\nabla^p R \cdot V_p - \frac{1}{2}\frac{\partial}{\partial t}|V|^2 + R_{pq}V^p V^q,\end{aligned}$$

where to obtain the last equality we used

$$-g^{pq}\frac{\partial}{\partial t}V_q \cdot V_p = -\frac{1}{2}\frac{\partial}{\partial t}|V|^2 + (R_{pq} - \nabla_p V_q)V^p V^q.$$

Combining the two equations above yields

$$\tilde{R}_{00} - \tilde{\nabla}_0 \tilde{V}_0 = \frac{1}{2}\frac{\partial}{\partial t}R - \frac{\partial}{\partial t}h + \frac{1}{2}\frac{\partial}{\partial t}|V|^2,$$

which is the same as formula 2. The proof of the lemma is complete.

## 4  Linearizing the Ricci flow and the space-time connection

In this section, we linearize the Ricci flow. In doing so, we must be careful to apply DeTurck's trick (see [D]) of modifying by a Lie derivative of the metric term (action of an infinitesimal diffeomorphism on the metric) in the variation direction only. We then have the following properties relevant to the linear trace Harnack inequality. The variation of the metric satisfies (1.3) (see Lemma 4.1, (ii) below.) The equation for the variation of the space-time Christoffel symbols $\tilde{\Gamma}^k_{ij}$ (defined in the last section) is formally the same as the equation for the variation of the Levi-Civita connection, where $g, \nabla$, and $h$ are replaced by $\tilde{g}, \tilde{\nabla}$, and $\tilde{h}$ (see Theorem 4.2 below.) In the next section, using this setup, we explain why the linear trace Harnack quadratic has such a nice evolution equation.

Let $g_{ij}(t)$ be a solution to the Ricci flow

$$\frac{\partial}{\partial t}g_{ij} = -2R_{ij}$$
$$g_{ij}(0) = (g_0)_{ij}$$

and $h_{ij}(t)$ be a solution to

$$\frac{\partial}{\partial t}h_{ij} = \Delta h_{ij} + 2R_{pijq}h_{pq} - R_{iq}h_{jq} - R_{jq}h_{iq}, \qquad (4.1)$$
$$h_{ij}(0) = (h_0)_{ij}.$$

The tensors $h(t)$ may be considered as the variation of the metrics $g(t)$ by a family of solutions to the modified Ricci flow. In particular, consider a 2-parameter family of metrics $g_{ij}(t,s)$ such that

$$\frac{\partial}{\partial t}g_{ij}(t,s) = (-2R_{ij} + \nabla_i W_j + \nabla_j W_i)(t,s)$$
$$g(0,0) = g_0,$$

where
$$W^k(t,s) = g^{ij}(t,s)(\Gamma_{ij}^k(t,s) - \Gamma_{ij}^k(t,0)) \qquad (4.2)$$

(here $\Gamma_{ij}^k(t,s)$ denote the Christoffel symbols of the metric $g_{ij}(t,s)$,) and

$$\frac{\partial}{\partial s}g_{ij}(0,0) = (h_0)_{ij}.$$

**Lemma 4.1**

(i) $g_{ij}(t,0) = g_{ij}(t)$

(ii) $h_{ij}(t,0) = \frac{\partial}{\partial s}g_{ij}(t,0)$.

PROOF. (i). This follows from the fact $W^k(t,0) = 0$, $g(0,0) = g_0$, and the uniqueness of solutions to the Ricci flow.

(ii). We compute

$$\frac{\partial}{\partial t}(\frac{\partial}{\partial s}g_{ij})(t,0) = \frac{\partial}{\partial s}(\frac{\partial}{\partial t}g_{ij})(t,0)$$
$$= \frac{\partial}{\partial s}(-2R_{ij} + \nabla_i W_j + \nabla_j W_i)(t,0)$$
$$= \left[\Delta(\frac{\partial}{\partial s}g_{ij}) + 2R_{kijl} \cdot \frac{\partial}{\partial s}g_{kl} - R_{ik} \cdot \frac{\partial}{\partial s}g_{kj} - R_{jk} \cdot \frac{\partial}{\partial s}g_{ki}\right](t,0).$$

Since both $\frac{\partial}{\partial s}g_{ij}(t,0)$ and $h_{ij}(t)$ satisfy the same equation (4.1) and $\frac{\partial}{\partial s}g_{ij}(0,0) = h_{ij}(0)$, by the uniqueness of solutions to (4.1) with given initial condition, we have

$$\frac{\partial}{\partial s}g_{ij}(t,0) = h_{ij}(t),$$

and the lemma is proved.

Since the family $\{g(t,s)\}_{t\in[0,T)}$ for $s$ fixed is a solution to the modified Ricci flow, we may associate to it the degenerate metric $\tilde{g}(t,s)$ and connection $\tilde{\nabla}(t,s)$ on space-time $M \times [0,T) \times \{t\}$ defined in section 3,

$$\tilde{\Gamma}^k_{ij} = \Gamma^k_{ij} \text{ if } i,j,k \geq 1$$
$$\tilde{\Gamma}^0_{ij} = 0 \text{ if } i,j \geq 0$$
$$\tilde{\Gamma}^k_{i0} = -R^k_i + \nabla_i W^k \text{ if } i,k \geq 1$$
$$\tilde{\Gamma}^k_{00} = -\tfrac{1}{2}\nabla^k(R+|W|^2) + g^{kp}\tfrac{\partial}{\partial t}W_p \text{ if } k \geq 1,$$

where $W$ is given by (4.2).

The tensor $\tilde{h}$ may be considered as the variation of the pair $(\tilde{g}, \tilde{\nabla})$ in the following sense.

**Theorem 4.2** *The space-time metric $\tilde{g}$ and connection $\tilde{\nabla}$ satisfy*

$$\frac{\partial}{\partial s}\tilde{g}^{ij} = -\tilde{g}^{ik}\tilde{g}^{jl}\tilde{h}_{kl}$$

*and*

$$\frac{\partial}{\partial s}\tilde{\Gamma}^k_{ij} = \frac{1}{2}\tilde{g}^{kl}(\tilde{\nabla}_i\tilde{h}_{jl} + \tilde{\nabla}_j\tilde{h}_{il} - \tilde{\nabla}_l\tilde{h}_{ij})$$

*at any point in $M \times [0,T) \times \{0\}$.*

PROOF. The variation of the metric is obvious, so we consider the variation of the connection. When $i,j,k \geq 1$, this follows from $\frac{\partial}{\partial s}g_{ij} = h_{ij}$ and the standard formula

$$\frac{\partial}{\partial s}\Gamma^k_{ij} = \frac{1}{2}g^{kl}\left(\nabla_i(\frac{\partial}{\partial s}g_{jl}) + \nabla_j(\frac{\partial}{\partial s}g_{il}) - \nabla_l(\frac{\partial}{\partial s}g_{ij})\right).$$

For $j = 0$ and $i,k \geq 1$, we compute

$$\frac{\partial}{\partial s}\tilde{\Gamma}^k_{i0} = \frac{\partial}{\partial s}\left[g^{kp}(-R_{ip} + \nabla_i W_p)\right]$$
$$= h_{kp}R_{ip} - \frac{\partial}{\partial s}R_{ik} + \nabla_i(\frac{\partial}{\partial s}W_k),$$

where we use the fact that $W(t,0) \equiv 0$.

Now at any point in $M \times [0,T) \times \{0\}$, we have

$$\frac{\partial}{\partial s}W_k = g_{kp}\frac{\partial}{\partial s}W^p = g_{kp}g^{ij}\frac{\partial}{\partial s}\Gamma^p_{ij} \qquad (4.3)$$
$$= \frac{1}{2}g^{ij}(\nabla_i h_{jk} + \nabla_j h_{ik} - \nabla_k h_{ij})$$
$$= div(h)_k - \frac{1}{2}\nabla_k H.$$

We also compute

$$\frac{\partial}{\partial s}R_{ik} = \frac{1}{2}(\nabla_q\nabla_i h_{kq} + \nabla_q\nabla_k h_{iq} - \Delta h_{ik} - \nabla_i\nabla_k H)$$
$$= \frac{1}{2}\nabla_i div(h)_k + \frac{1}{2}\nabla_k div(h)_i - R_{qikp}h_{pq} + \frac{1}{2}R_{ip}h_{kp} + \frac{1}{2}R_{kp}h_{ip}$$
$$- \frac{1}{2}\Delta h_{ik} - \nabla_i\nabla_k H.$$

Combining all of the above equations, we find that

$$\frac{\partial}{\partial s}\tilde{\Gamma}^k_{i0} = \frac{1}{2}\nabla_i div(h)_k - \frac{1}{2}\nabla_k div(h)_i + \frac{1}{2}\Delta h_{ik} + R_{pikq}h_{pq} + \frac{1}{2}R_{ip}h_{pk} - \frac{1}{2}R_{kp}h_{pi}.$$

Next we compute

$$\frac{1}{2}\tilde{g}^{kl}(\tilde{\nabla}_i\tilde{h}_{0l} + \tilde{\nabla}_0\tilde{h}_{il} - \tilde{\nabla}_l\tilde{h}_{i0})$$
$$= \frac{1}{2}[\nabla_i div(h)_k - \tilde{\Gamma}^p_{i0}h_{pk} + \frac{\partial}{\partial t}h_{ik} - \tilde{\Gamma}^p_{0i}h_{pk} - \tilde{\Gamma}^p_{0k}h_{ip} - \nabla_k div(h)_i + \tilde{\Gamma}^p_{k0}h_{ip}]$$
$$= \frac{1}{2}[\nabla_i div(h)_k - \nabla_k div(h)_i + \Delta h_{ik} + 2R_{pikq}h_{pq} + R_{ip}h_{pk} - R_{kp}h_{pi}]$$
$$= \frac{\partial}{\partial s}\tilde{\Gamma}^k_{i0}.$$

Finally, for $i = j = 0$ and $k \geq 1$ we compute

$$\frac{\partial}{\partial s}\tilde{\Gamma}^k_{00} = \frac{1}{2}h_{kp}\nabla_p R - \frac{1}{2}\nabla^k\left(\frac{\partial}{\partial s}R\right) + g^{kp}\frac{\partial}{\partial s}\left(\frac{\partial}{\partial t}W_p\right),$$

where we used $W(t,0) = 0$ and $\frac{\partial}{\partial s}W(t,0) = 0$.

Since

$$\frac{\partial}{\partial s}R = div(div(h)) - \Delta H - Rc \cdot h \tag{4.4}$$

and

$$\frac{\partial}{\partial s}\left(\frac{\partial}{\partial t}W_p\right) = \frac{\partial}{\partial t}\left(\frac{\partial}{\partial s}W_p\right) = \frac{\partial}{\partial t}(div(h)_p - \frac{1}{2}\nabla_p H),$$

we have

$$\frac{\partial}{\partial s}\tilde{\Gamma}^k_{00} = \frac{1}{2}h_{kp}\nabla_p R - \frac{1}{2}\nabla_k\left[div(div(h)) - \Delta H - Rc \cdot h\right] + \frac{\partial}{\partial t}(div(h)_k - \frac{1}{2}\nabla_k H).$$

Now from (4.1), we compute

$$\frac{\partial}{\partial t}H = \Delta H + 2Rc \cdot h, \tag{4.5}$$

which implies

$$\frac{\partial}{\partial s}\tilde{\Gamma}^k_{00} = \frac{1}{2}h_{kp}\nabla_p R - \frac{1}{2}\nabla_k\left[div(div(h)) + Rc \cdot h\right] + \frac{\partial}{\partial t}(div(h))_k.$$

On the other hand, we compute

$$\frac{1}{2}\tilde{g}^{kl}(\tilde{\nabla}_0\tilde{h}_{0l} + \tilde{\nabla}_0\tilde{h}_{0l} - \tilde{\nabla}_l\tilde{h}_{00})$$
$$= \frac{\partial}{\partial t}(div(h)_k) - \tilde{\Gamma}^p_{00}h_{pk} - \tilde{\Gamma}^p_{0k}\tilde{h}_{0p} - \frac{1}{2}\nabla_k\left[div(div(h)) + Rc \cdot h\right] + \tilde{\Gamma}^p_{k0}\tilde{h}_{0p}$$
$$= \frac{\partial}{\partial t}(div(h)_k) + \frac{1}{2}\nabla^p R h_{pk} - \frac{1}{2}\nabla_k\left[div(div(h)) + Rc \cdot h\right]$$
$$= \frac{\partial}{\partial s}\tilde{\Gamma}^k_{00}.$$

The proof of the lemma is complete.

## 5 Linear trace Harnack quadratic and linearizing the space-time Ricci tensor

In this section we linearize the space-time Ricci tensor $\tilde{R}_{ij}$, given a solution to the Ricci flow

$$\frac{\partial}{\partial t}g_{ij} = -2R_{ij},$$

and then consider the evolution equation for the linear trace Harnack quadratic from this point of view. We consider the same 2-parameter family of metrics $\{g(t,s)\}$ as before. Again, we use DeTurck's trick (now in space-time) and add a Lie derivative of the metric term to the Ricci tensor and consider

$$-2\tilde{R}_{ij} + \tilde{\nabla}_i\tilde{W}_j + \tilde{\nabla}_j\tilde{W}_i.$$

Clearly, we want $\tilde{W}$ to be an extension of $W$:

$$\tilde{W}^k(t,s) = W^k(t,s) = g^{ij}(t,s)(\Gamma^k_{ij}(t,s) - \Gamma^k_{ij}(t,0)),$$

for $k \geq 1$.

The key is to define $\tilde{W}_0$ suitably, unlike in Theorem 3.1 for $\tilde{V}_0$, we cannot prescribe $\tilde{W}_0$ arbitrarily. The right definition is to let

$$\tilde{W}_0 = \sum_{p=1}^{n}\nabla^p W_p,$$

or equivalently [1]
$$\tilde{W}_0 = \tilde{g}^{pq}\tilde{\nabla}_p \tilde{W}_q. \tag{5.1}$$

We then have

**Theorem 5.1** *On $M \times [0,T) \times \{0\}$*

$$\frac{\partial}{\partial s}(-2\tilde{R}_{ij} + \tilde{\nabla}_i \tilde{W}_j + \tilde{\nabla}_j \tilde{W}_i) = \tilde{\Delta}\tilde{h}_{ij} + 2\tilde{R}^q_{pij}\tilde{h}_{pq} - \tilde{R}_{iq}\tilde{h}_{jq} - \tilde{R}_{jq}\tilde{h}_{iq}.$$

PROOF. We first compute

$$\frac{\partial}{\partial s}\tilde{R}_{ij} = \tilde{\nabla}_p(\frac{\partial}{\partial s}\tilde{\Gamma}^p_{ij}) - \tilde{\nabla}_i(\frac{\partial}{\partial s}\tilde{\Gamma}^p_{pj})$$

$$= \frac{1}{2}\tilde{\nabla}_p(\tilde{\nabla}_i\tilde{h}_{jp} + \tilde{\nabla}_j\tilde{h}_{ip} - \tilde{\nabla}_p\tilde{h}_{ij}) - \frac{1}{2}\tilde{\nabla}_i(\tilde{\nabla}_p\tilde{h}_{jp} + \tilde{\nabla}_j\tilde{h}_{pp} - \tilde{\nabla}_p\tilde{h}_{pj})$$

$$= -\frac{1}{2}\tilde{\Delta}\tilde{h}_{ij} - \tilde{R}^q_{pij}\tilde{h}_{pq} + \frac{1}{2}\tilde{R}_{iq}\tilde{h}_{jq} + \frac{1}{2}\tilde{R}_{jq}\tilde{h}_{iq} + \frac{1}{2}\tilde{\nabla}_i\tilde{\nabla}_p\tilde{h}_{jp} + \frac{1}{2}\tilde{\nabla}_j\tilde{\nabla}_p\tilde{h}_{ip}$$

$$- \frac{1}{2}\tilde{\nabla}_i\tilde{\nabla}_j\tilde{h}_{pp}.$$

Hence the theorem will follow from showing that

$$\frac{\partial}{\partial s}\tilde{W}_j = \tilde{\nabla}_p\tilde{h}_{jp} - \frac{1}{2}\tilde{\nabla}_j\tilde{h}_{pp},$$

for $j \geq 0$. For $j \geq 1$, this is true by equation (4.3). For $j = 0$, this follows from the computation

$$\frac{\partial}{\partial s}\tilde{W}_0 = \nabla^p(\frac{\partial}{\partial s}W^p) = div(div(h)) - \frac{1}{2}\Delta H \tag{5.2}$$

$$= \tilde{\nabla}_p\tilde{h}_{0p} - \frac{1}{2}\tilde{\nabla}_0\tilde{h}_{pp},$$

where we used definition (5.1), $W = 0$ on $M \times [0,T) \times \{0\}$, and (4.5). The theorem is proved.

This theorem together with Lemma 2.1, says that

$$\frac{\partial}{\partial t}\tilde{h}_{ij} = \tilde{\Delta}_L \tilde{h}_{ij} = \frac{\partial}{\partial s}(-2\tilde{R}_{ij} + \tilde{\nabla}_i\tilde{W}_j + \tilde{\nabla}_j\tilde{W}_i), \tag{5.3}$$

---

[1] It is perhaps natural to propose that all covariant space tensors be extended this way. In particular, extending the Riemann curvature and Ricci tensors of $g$ this way yield the Riemann curvature and Ricci tensors of the space-time connection $\tilde{\nabla}$. In addition, the extension of $h$ to $\tilde{h}$ is given in exactly this way.

where by definition,
$$\tilde{\Delta}_L \tilde{h}_{ij} = \tilde{\Delta} \tilde{h}_{ij} + 2\tilde{R}^q_{pij} \tilde{h}_{pq} - \tilde{R}_{iq} \tilde{h}_{jq} - \tilde{R}_{jq} \tilde{h}_{iq}.$$

This formula further confirms that $\tilde{h}$ should be thought of as the variation of the pair $(\tilde{g}, \tilde{\nabla})$. It also agrees with the equation obtained by comparing the mixed partial derivatives of the space-time connection in the time and variation direction. In particular, first recall that

$$\frac{\partial}{\partial s} \tilde{\Gamma}^k_{ij} = \frac{1}{2} \tilde{g}^{kl} (\tilde{\nabla}_i \tilde{h}_{jl} + \tilde{\nabla}_j \tilde{h}_{il} - \tilde{\nabla}_l \tilde{h}_{ij}).$$

Taking $k \geq 1$, we compute (omitting details:)

$$\frac{\partial}{\partial t}(\frac{\partial}{\partial s} \tilde{\Gamma}^k_{ij}) = \sum_{l \geq 1} \tilde{R}_{kl}(\tilde{\nabla}_i \tilde{h}_{jl} + \tilde{\nabla}_j \tilde{h}_{il} - \tilde{\nabla}_l \tilde{h}_{ij}) + \sum_{p \geq 1} (\tilde{\nabla}_i \tilde{R}_{jp} + \tilde{\nabla}_j \tilde{R}_{ip} - \tilde{\nabla}_p \tilde{R}_{ij}) \tilde{h}_{pk}$$

$$+ \frac{1}{2} \left[ \tilde{\nabla}_i (\frac{\partial}{\partial t} \tilde{h}_{jk}) + \tilde{\nabla}_j (\frac{\partial}{\partial t} \tilde{h}_{ik}) - \tilde{\nabla}_k (\frac{\partial}{\partial t} \tilde{h}_{ij}) \right]. \tag{5.4}$$

Next, we compute the same quantity, except that we switch the order of the partial derivatives. We have (again omitting the details of the computation:)

$$\frac{\partial}{\partial s}(\frac{\partial}{\partial t} \tilde{\Gamma}^k_{ij}) = \sum_{l \geq 1} \tilde{h}_{kl}(\tilde{\nabla}_i \tilde{R}_{jl} + \tilde{\nabla}_j \tilde{R}_{il} - \tilde{\nabla}_l \tilde{R}_{ij}) + \sum_{p \geq 1} (\tilde{\nabla}_i \tilde{h}_{jp} + \tilde{\nabla}_j \tilde{h}_{ip} - \tilde{\nabla}_p \tilde{h}_{ij}) \tilde{R}_{pk}$$

$$+ \frac{1}{2} \left[ \tilde{\nabla}_i \tilde{S}_{jk} + \tilde{\nabla}_j \tilde{S}_{ik} - \tilde{\nabla}_k \tilde{S}_{ij} \right]. \tag{5.5}$$

where
$$\tilde{S}_{ij} = \frac{\partial}{\partial s} \left( -2\tilde{R}_{ij} + \tilde{\nabla}_i \tilde{W}_j + \tilde{\nabla}_j \tilde{W}_i \right) = \tilde{\Delta}_L \tilde{h}_{ij}.$$

Comparing (5.4) and (5.5), we obtain

$$\tilde{\nabla}_i (\frac{\partial}{\partial t} \tilde{h}_{jk} - \tilde{S}_{jk}) + \tilde{\nabla}_j (\frac{\partial}{\partial t} \tilde{h}_{ik} - \tilde{S}_{ik}) - \tilde{\nabla}_k (\frac{\partial}{\partial t} \tilde{h}_{ij} - \tilde{S}_{ij}) = 0, \tag{5.6}$$

for all $i, j \geq 0$ and $k \geq 1$. This agrees with (5.3); note however, (5.6) doesn't imply (2.5). To further understand the interpretation of $\tilde{h}$ as a variation of $(\tilde{g}, \tilde{\nabla})$, we consider the approximation approach of the next section.

# 6  An approximation approach

In [CC], the authors showed that there is a family of Riemannian metrics $\bar{g}$ on $\tilde{M} = M \times [0, T)$ such that $\bar{g}^{-1}$ limits to the degenerate metric $\tilde{g}$ on $T^*\tilde{M}$ and the Levi-Civita connection $\bar{\Gamma}_{ij}^k$ limits to the space-time connection $\tilde{\Gamma}_{ij}^k$. Here, we adopt this approach to show that, corresponding to the 2-parameter family of metrics $\{g(t,s)\}$ defined in section 4, there is a family of Riemannian metrics $\widehat{g}$, parametrized by $s$ and a sufficiently large number $N$, on $\tilde{M}$ which satisfy the above properties, and also the additional property that the variation in the $s$ direction of each metric $\widehat{g}$ in the family is the linear trace Harnack quantity $\tilde{h}$. Furthermore, we derive the evolution equation for $\tilde{h}$ again using these approximate metrics.

In general, given a solution $g(t)$ to the modified Ricci flow

$$\frac{\partial}{\partial t} g_{ij} = -2 R_{ij} + \nabla_i V_j + \nabla_j V_i,$$

and a function $f : M \times [0, T) \to \mathbb{R}$, define the Riemannian metric $\widehat{g}$ (parametrized by large numbers $N$) on $\tilde{M} = M \times [0, T)$ by

$$\begin{aligned}
\widehat{g}_{ij} &= g_{ij} & \text{if } i, j \geq 1 \\
\widehat{g}_{i0} &= V_i + \nabla_i f & \text{if } i \geq 1 \\
\widehat{g}_{00} &= R + |V|^2 + 2\tfrac{\partial f}{\partial t} + N,
\end{aligned}$$

where $N$ is a sufficiently large positive real number to make $\widehat{g}$ positive-definite. [2] Clearly

$$\lim_{N \to \infty} \widehat{g}^{-1} = \tilde{g}.$$

We shall show that for any function $f$, the Levi-Civita connections of this family of Riemannian metrics converge (as $N \to \infty$) to the space-time connection $\tilde{\Gamma}_{ij}^k$. We then apply this to the 2-parameter family of metrics defined in section 4. First recall that if $I$ denotes the $n \times n$ identity matrix, $X$ is a column $n$-vector, and $b$ is real number greater than $|X|^2$, then

$$\begin{pmatrix} I & X \\ X^T & b \end{pmatrix}^{-1} = \begin{pmatrix} I + \frac{X X^T}{b - |X|^2} & -\frac{1}{b - |X|^2} X \\ -\frac{1}{b - |X|^2} X^T & \frac{1}{b - |X|^2} \end{pmatrix}.$$

---

[2] Conan Leung pointed out to us that this construction of a family of positive-definite metrics (together with their Levi-Civita connections) limiting to a metric which is infinite in certain directions is similar to a construction in work of Bismut [B].

Hence the inverse of the metric $\widehat{g}$ is given by

$$\widehat{g}^{ij} = g^{ij} + \frac{(V^i+\nabla^i f)(V^j+\nabla^j f)}{R+2\left(\frac{\partial f}{\partial t}-\langle\nabla f,V\rangle\right)-|\nabla f|^2+N} \quad \text{if } i,j \geq 1$$

$$\widehat{g}^{i0} = -\frac{1}{R+2\left(\frac{\partial f}{\partial t}-\langle\nabla f,V\rangle\right)-|\nabla f|^2+N}\left(V^i+\nabla^i f\right) \quad \text{if } i \geq 1$$

$$\widehat{g}^{00} = \frac{1}{R+2\left(\frac{\partial f}{\partial t}-\langle\nabla f,V\rangle\right)-|\nabla f|^2+N}.$$

Using the standard formula

$$\widehat{\Gamma}^k_{ij} = \frac{1}{2}\sum_{l=0}^n \widehat{g}^{kl}\left(\partial_i \widehat{g}_{jl} + \partial_j \widehat{g}_{il} - \partial_l \widehat{g}_{ij}\right),$$

we compute that the Christoffel symbols are given by (the details of the proof are omitted:)

**Lemma 6.1** *For $i,j,k \geq 1$, we have*

$$\widehat{\Gamma}^k_{ij} = \Gamma^k_{ij} - \frac{1}{R+2\left(\frac{\partial f}{\partial t}-\langle\nabla f,V\rangle\right)-|\nabla f|^2+N}\left(V^k+\nabla^k f\right)(R_{ij}+\nabla_i\nabla_j f)$$

$$\widehat{\Gamma}^k_{i0} = -R^k_i + \nabla_i V^k - \frac{(V^k+\nabla^k f)\left(\frac{1}{2}\nabla_i R + V^l R_{il} + \nabla_i \frac{\partial f}{\partial t} - \nabla^l f(-R_{il}+\nabla_i V_l)\right)}{R+2\left(\frac{\partial f}{\partial t}-\langle\nabla f,V\rangle\right)-|\nabla f|^2+N}$$

$$\widehat{\Gamma}^k_{00} = g^{kl}\frac{\partial}{\partial t}V_l - \frac{1}{2}\nabla^k\left(R+|V|^2\right)$$
$$- \frac{(V^k+\nabla^k f)}{R+2\left(\frac{\partial f}{\partial t}-\langle\nabla f,V\rangle\right)-|\nabla f|^2+N}\left[\begin{array}{c}\frac{1}{2}\sum_{l=1}^n V^l\nabla_l R + \frac{1}{2}\frac{\partial R}{\partial t} + \sum_{i,j=1}^n R_{ij}V_iV_j \\ -\nabla^l f\left(\frac{\partial}{\partial t}V_l - \frac{1}{2}\partial_l\left(R+|V|^2\right)\right) + \frac{\partial}{\partial t}\left(\frac{\partial f}{\partial t}\right)\end{array}\right]$$

$$\widehat{\Gamma}^0_{ij} = \frac{1}{R+2\left(\frac{\partial f}{\partial t}-\langle\nabla f,V\rangle\right)-|\nabla f|^2+N}(R_{ij}+\nabla_i\nabla_j f)$$

$$\widehat{\Gamma}^0_{i0} = \frac{\left[\frac{1}{2}\nabla_i R + \sum_{l=1}^n V^l R_{il} + \sum_{l=1}^n \nabla^l f(R_{il}-\nabla_i V_l) + \nabla_i \frac{\partial f}{\partial t}\right]}{R+2\left(\frac{\partial f}{\partial t}-\langle\nabla f,V\rangle\right)-|\nabla f|^2+N}$$

$$\widehat{\Gamma}^0_{00} = \frac{\left(\frac{1}{2}\sum_{l=1}^n V^l\nabla_l R + \frac{1}{2}\frac{\partial}{\partial t}R + R_{ij}V_iV_j - \nabla^l f\left(\frac{\partial}{\partial t}(V_l+\nabla_l f) - \frac{1}{2}\nabla_l\left(R+|V|^2+2\frac{\partial f}{\partial t}\right)\right)\right)}{R+2\left(\frac{\partial f}{\partial t}-\langle\nabla f,V\rangle\right)-|\nabla f|^2+N}$$

It follows that

**Lemma 6.2** *For $N$ sufficiently large*

$$\widehat{g}^{ij} = \tilde{g}^{ij} + O\left(\frac{1}{N}\right)$$

$$\widehat{\Gamma}^k_{ij} = \tilde{\Gamma}^k_{ij} + O\left(\frac{1}{N}\right)$$

$$\widehat{R}^l_{ijk} = \tilde{R}^l_{ijk} + O\left(\frac{1}{N}\right),$$

*for all $i, j, k \geq 0$. In particular,*

$$\lim_{N \to \infty} \widehat{g}^{ij} = \tilde{g}^{ij}$$

$$\lim_{N \to \infty} \widehat{\Gamma}^k_{ij} = \tilde{\Gamma}^k_{ij}$$

$$\lim_{N \to \infty} \widehat{R}^l_{ijk} = \tilde{R}^l_{ijk}.$$

Now let $\{g(t,s)\}$ denote the 2-parameter family of metrics defined in section 4, which satisfy the modified Ricci flow in the $t$ direction

$$\frac{\partial}{\partial t} g_{ij} = -2R_{ij} + \nabla_i W_j + \nabla_j W_i,$$

where $W$ is defined by (4.2). Furthermore, take

$$f = \frac{1}{2} \log \frac{\det g(t,s)}{\det g(t,0)} = \log \frac{dV(t,s)}{dV(t,0)},$$

so that

$$\widehat{g}_{ij} = g_{ij}$$

$$\widehat{g}_{i0} = W_i + \frac{1}{2} \nabla_i \log \frac{\det g(t,s)}{\det g(t,0)}$$

$$\widehat{g}_{00} = R + |W|^2 + \frac{\partial}{\partial t} \log \frac{\det g(t,s)}{\det g(t,0)} + N.$$

In addition to the properties of Lemma 6.2, we now have that the variation of $\widehat{g}$ is the linear trace Harnack quantity $\tilde{h}$

**Lemma 6.3** *For $i, j \geq 0$,*

$$\frac{\partial}{\partial s} \widehat{g}_{ij} = \tilde{h}_{ij}.$$

PROOF. We have for $i, j \geq 1$, by definition
$$\frac{\partial}{\partial s}\widehat{g}_{ij} = \frac{\partial}{\partial s}g_{ij} = h_{ij} = \tilde{h}_{ij}.$$

For $i \geq 1$, we compute using equation (4.3) that
$$\frac{\partial}{\partial s}\widehat{g}_{i0} = \frac{\partial W_i}{\partial s} + \frac{1}{2}\nabla_i H = div(h)_i = \tilde{h}_{i0},$$

and using equations (4.4) and (4.5) that
$$\frac{\partial}{\partial s}\widehat{g}_{00} = \frac{\partial R}{\partial s} + \frac{\partial H}{\partial t} = div\,(div(h)) + Rc \cdot h = \tilde{h}_{00}.$$

We also have, to first order in $s$, that the time-derivative of $\widehat{g}$ is the modified space-time Ricci tensor.

**Lemma 6.4**
$$\frac{\partial}{\partial t}\widehat{g}_{ij} = -2\tilde{R}_{ij} + \tilde{\nabla}_i\tilde{W}_j + \tilde{\nabla}_j\tilde{W}_i + o(s),$$

*where*
$$\left.\frac{\partial}{\partial s}\right|_{s=0} o(s) = 0.$$

Using Theorem 5.1 and Lemmas 6.3 and 6.4, we can give the following alternate proof of Lemma 2.1 (or equivalently Proposition 2.3.) On $M \times [0, T) \times \{0\}$

$$\frac{\partial}{\partial t}\tilde{h}_{ij} = \frac{\partial^2}{\partial t \partial s}\widehat{g}_{ij} = \frac{\partial^2}{\partial s \partial t}\widehat{g}_{ij} = \left.\frac{\partial}{\partial s}\right|_{s=0}\left(-2\tilde{R}_{ij} + \tilde{\nabla}_i\tilde{W}_j + \tilde{\nabla}_j\tilde{W}_i\right) = \tilde{\triangle}_L \tilde{h}_{ij}.$$

**Acknowledgement.** We would like to thank Richard Hamilton for very helpful discussions.

# References

[B] J.M. Bismut, The index theorem for families of Dirac operators: two heat equation proofs, Invent. Math. **83**. p.91-151 (1986).

[CC] B. Chow and S.C. Chu, A geometric interpretation of Hamilton's Harnack inequality for the Ricci flow, Math. Research Letters **2**. p.701-718 (1995).